\documentclass[reqno,12pt,a4]{amsart}
\usepackage{amssymb}
\usepackage{amsmath}

\usepackage[all]{xy}

\def\Ext{\mathop{\rm Ext}\nolimits}
\def\id{\mathop{\rm id}\nolimits}

\def\b{{\beta}}

\def\s{{\sigma}}

\newtheorem{thm}{Theorem}[section]
\newtheorem{lem}[thm]{Lemma}
\newtheorem{prop}[thm]{Proposition}
\newtheorem{dfn}[thm]{Definition}
\newtheorem{cor}[thm]{Corollary}

\def\K{{\mathcal K}}

\date{}
\author{V. Manuilov and K. Thomsen}
\title{On the asymptotic tensor norm}

\begin{document}

 \begin{abstract}
We introduce a new asymptotic one-sided and symmetric tensor norm,
the latter of which can be considered as the minimal tensor norm on the
category of separable $C^*$-algebras with homotopy classes of asymptotic
homomorphisms as morphisms. We show that the one-sided asymptotic tensor
norm differs in general from both the minimal and the maximal tensor norms
and discuss its relation to semi-invertibility of $C^*$-extensions.

 \end{abstract}

\maketitle

\section{Introduction. Invertibility of extensions}

One of the reasons, why the Brown--Douglas--Fillmore theory~\cite{BDF}
gives so nice a classification for extensions of nuclear $C^*$-algebras is
their invertibility. Beyond the nuclear case not much is known about
general classification of $C^*$-extensions, but more and more examples of
non-invertible extensions are coming up~\cite{Anderson,Kirchberg,
Wass1,Wass2,Haagerup}. In~\cite{MT2,MT3} it was suggested to weaken the
notion of triviality for extensions so that more extensions would become
invertible in this new sense. As it was shown by S.~Wassermann, one of the
reasons for non-invertibility is non-exactness and his idea relates
many of examples of non-invertible extensions to the problem of
coincidence for certain tensor product norms. In the present paper we
develop this idea to define the asymptotic tensor norm and study this
norm in hope to learn more about extensions.

Let
\begin{equation}\label{extension}
\begin{xymatrix}{
0\ar[r]& B\ar[r]& E\ar[r]^-{p}& A\ar[r]& 0
}
\end{xymatrix}
\end{equation}
be an extension of $A$ by $B$, i.e. a short
exact sequence of $C^*$-algebras. We always assume that $B$ is stable,
i.e. $B\otimes\K\cong B$, where $\K$ is the $C^*$-algebra of compact
operators. We aslo assume all $C^*$-algebras to be separable with some
obvious exceptions like the algebra $L(H)$ of bounded operators on a
Hilbert space $H$. The extension (\ref{extension}) is called {\it split}
if there is a $*$-homomorphism $s:A\to E$ that is the right inverse for
the surjection $p$, i.e. $p\circ s=\id_A$. Among properties of split
extensions one should mention their homotopy triviality (i.e. homotopy
triviality of their Busby invariant). Due to stability of $B$, one can fix
an isomorphism $M_2(B)\cong B$, which makes it possible to define a direct
sum of two extensions of $A$ by $B$. An extension of $A$ by $B$ is called
{\it invertible} if there is another extension of $A$ by $B$, such that
their direct sum is split. Remark that, due to the Stinespring theorem,
this is equivalent to existence of a completely positive splitting for $p$.
Wassermann's idea of producing non-invertible extensions works as follows.
Take a non-exact $C^*$-algebra $D$ and an extension of the form
(\ref{extension}) such that the sequence
 $$
0\to B\otimes_{\min}D\to E\otimes_{\min}D\to A\otimes_{\min}D\to 0
 $$
is not exact and denote by $A\otimes_E D$ the completion of the algebraic
tensor product $A\odot D$ given by the quotient
$E\otimes_{\min}D/B\otimes_{\min}D$. By $\|\cdot\|_E$ we denote the
corresponding tensor norm on $A\otimes_E D$.

\begin{thm}[S.~Wassermann,~\cite{Wass-book}]
If the norm $\|\cdot\|_E$ on $A\odot D$ differs from the minimal tensor
norm then the extension $(\ref{extension})$ is not invertible.

\end{thm}

Remark that for any commutative $C^*$-algebra $C_0(X)$ one can form an
extension
 $$
0\to C_0(X)\otimes B\to C_0(X)\otimes E\to C_0(X)\otimes A\to 0
 $$
out of the extension (\ref{extension}).
It is easy to see that, if $\|\cdot\|_E\neq\|\cdot\|_{\min}$ on $A\odot D$,
then $\|\cdot\|_{C_0(X)\otimes E}\neq \|\cdot\|_{\min}$ on $C_0(X)\otimes
A\odot D$. Taking $X=[0,1)$, one obtains plenty of examples of
non-invertible extensions of {\it contractible} $C^*$-algebras. The first
such example was found in~\cite{Kirchberg} (in fact,
Kirchberg's example is much more interesting than these ones because there 
one has $B=\K$).

In order to have more invertible extensions one has to change definitions.
Such a change was suggested in~\cite{MT2,MT3}:
\begin{dfn}
{\rm
An extension (\ref{extension}) is {\it asymptotically split} if there
exists an asymptotic homomorphism $\s_t:A\to E$ that is the right inverse for
the surjection $p$ for every $t$, i.e. $p\circ s_t=\id_A$. An extension
(\ref{extension}) is {\it semi-invertible} if there is an extension $0\to
B\to E'\to A\to 0$, such that their direct sum is asymptotically split.

}
\end{dfn}

It is still easy to see that asymptotically split extensions are homotopy
trivial. But here we also have an almost inverse statement:

\begin{thm}[\cite{MT3}]\label{contractible}
Any extension of a contractible $C^*$-algebra is asymptotically split.

\end{thm}

Thus some examples of non-invertible extensions turn out to be
semi-invertible.

\section{Asymptotic tensor norms}

Let $L(H)$ denote the algebra of bounded operators on a
separable Hilbert space $H$.
By an {\it asymptotic representation} of a $C^*$-algebra $A$ on $H$
we mean an asymptotic homomorphism
$\mu = (\mu_t)_{t \in [0,\infty)} : A \to L(H)$. Note that a
genuine representation of $A$ on $H$ can be considered as an asymptotic
representation of $A$ in the obvious way. By using asymptotic
representations instead of genuine representations we shall now introduce
two new tensor norms on the algebraic tensor product $A \odot D$ of two
$C^*$-algebras $A$ and $D$.

Let $H_1,H_2$ be separable
Hilbert spaces and $\mu = (\mu_t)_{t \in [0,\infty)} :
A \to L(H_1), \nu = (\nu_t)_{t \in [0,\infty)} :
D \to L(H_2)$ two asymptotic representations. For each
$a \in A$ and $d \in D$, we can define elements $a^{\mu \odot \nu} ,
d^{\mu \odot \nu} \in C_b\left( [0,\infty),
L(H_1 \otimes H_2)\right)$ by
$$
a^{\mu \odot \nu}(t) = \mu_t(a) \otimes 1_{H_2}
\quad {\rm and}\quad
d^{\mu \odot \nu}(t) = 1_{H_1} \otimes \nu_t(d).
$$
We can then define a $*$-homomorphism
$$
\mu \odot \nu : A \odot D \to C_b\left( [0,\infty),
L(H_1 \otimes H_2)\right)/  C_0\left( [0,\infty),
L(H_1 \otimes H_2)\right)
$$
such that
$$
\mu \odot \nu\Bigl( \sum_i a_i \odot d_i\Bigr) =
\sum_i a_i^{\mu \odot \nu}d_i^{\mu \odot \nu} .
$$
We can then define a norm $\| \cdot \|_{\sigma}$ (the {\it symmetric
asymptotic tensor norm}) on $A \odot D$ by
$$
\Bigl\| \sum_i a_i \odot d_i \Bigr\|_{\sigma} = \sup_{\mu, \nu} 
\Bigl\|\mu \odot\nu\Bigl( \sum_i a_i \odot d_i\Bigr)\Bigr\|,
$$
where we take the supremum over all pairs $(\mu,\nu)$, where $\mu$ and
$\nu$ are asymptotic representations of $A$ and $D$, respectively. We
define also a norm $\| \cdot \|_{\lambda}$ (the {\it left asymptotic
tensor norm}) on $A \odot D$ by
$$
\Bigl\| \sum_i a_i \odot d_i \Bigr\|_{\lambda} = \sup_{\mu, \nu} 
\Bigl\|\mu \odot\nu\Bigl( \sum_i a_i \odot d_i\Bigr)\Bigr\|,
$$
where we take the supremum over all pairs $(\mu,\nu)$, where $\mu$ is an
asyptotic representation of $A$ and $\nu$ is a genuine representation of
$D$. Clearly, both asymptotic tensor norms, $\|\cdot\|_\lambda$ and 
$\|\cdot\|_\sigma$ are cross-norms and
$$
\| \cdot \|_{min} \leq \|\cdot \|_{\lambda} \leq \|\cdot \|_{\sigma}
\leq \|\cdot \|_{max} .
$$
We denote by $A \otimes_{\lambda} D$ and $A \otimes_{\sigma} D$ the
$C^*$-algebras obtained by completing $A \odot D$ with respect to the
norm $\| \cdot \|_{\lambda}$ and $\|\cdot \|_{\sigma}$ respectively.

\begin{lem}
Given $C^*$-algebras $A$ and $D$, there is an asymptotic representation
$\mu$ on a {\rm non-separable Hilbert space} of $A$, such that 
$\|\cdot\|_{\lambda}=\|\mu\odot\phi(\cdot)\|$, where $\phi$ is a
universal representation of $D$.

\end{lem}
\begin{proof}
Since any asymptotic representation is asymptotically equivalent to an
equicontinuous one, so in the definition of the norm $\|\cdot\|_\lambda$
one can use only equicontinuous asymptotic representations. 
Let $\{\pi^\alpha\}_{\alpha\in{\mathcal A}}$ be the set of all
equicontinuous asymptotic representations of $A$. We would like to take
direct sums, but this can be possible only after reparametrization. 
Notice that the norm $\|\nu\odot\phi(\sum_i a_i\odot d_i)\|$ does not
change after a reparametrization of an asymptotic homomorphism $\nu$.
Let $F_n\subset A$ be a sequence of finite subsets with dense union.
Fix a decreasing sequence $\{\varepsilon_n\}_{n\in{\mathbb N}}$ of
positive numbers vanishing at infinity. For each $\alpha$ let
$r_\alpha(t)$ be a reparametrization such that
$\|\pi^\alpha_{r_\alpha(t)}(ab)-\pi^\alpha_{r_\alpha(t)}(a)
\pi^\alpha_{r_\alpha(t)}(b)\|<\varepsilon_n$ for all $t\geq n$ when
$a,b\in F_n$ and similar conditions for other asymptotically algebraic
relations hold. Then $\mu_t=\oplus_{\alpha\in{\mathcal
A}}\pi^\alpha_{r_\alpha(t)}$ is an asymptotic representation of $A$ and
can be used in calculation of the norm $\|\cdot\|_\lambda$. 
\end{proof}

\begin{lem}
Let $\phi=(\phi_t)_{t\in[0,\infty)}:A_1\to A_2$ and
$\psi=(\psi_t)_{t\in[0,\infty)}:D_1\to D_2$ be asymptotic homomorphisms.
Then their tensor
product $\phi_t\otimes\psi_t$ extends to an asymptotic homomorphism from
$A_1\otimes_\sigma D_1$ to $A_2\otimes_\sigma D_2$.
If $\psi:D_1\to D_2$ is a genuine $*$-homomorphism then the tensor
product $\phi_t\otimes\psi$ extends to an asymptotic homomorphism from 
$A_1\otimes_\lambda D_1$ to $A_2\otimes_\lambda D_2$.
\end{lem}
\begin{proof}
Obvious.
\end{proof}

\begin{prop}\label{exact_seq}
For an extension (\ref{extension}),
suppose that there exists a $C^*$-algebra $D$ and an element $x\in A\odot
D$ such that $\|x\|_E >\|x\|_\lambda$. Then the extension
is not semi-invertible.

\end{prop}
\begin{proof}
The idea of the proof is borrowed from \cite{Wass2}.
Suppose the contrary, i.e. that (\ref{extension}) is semi-invertible. Then
there exists an extension 
$
\begin{xymatrix}{
0\ar[r]& B\ar[r]& E'\ar[r]^-{p'}& A\ar[r]& 0
}
\end{xymatrix}
$
and an asymptotic
splitting $s=(s_t)_{t\in{[0,\infty)}}:A\to C$, where $C\subset M_2(M(B))$
is the $C^*$-subalgebra of the form
$$
C=\left\lbrace\left(\begin{array}{cc}e&b_1\\b_2&e'\end{array}\right):
b_1,b_2\in B, e\in E,e'\in E', p(e)=p'(e')\right\rbrace. 
$$
Then we have a
well-defined asymptotic homomorphism
$s_t\otimes_\lambda \id_D:A\otimes_\lambda D\to C\otimes_{\min} D$.
Let $d:C\to E$ be a completely positive contraction given by 
$d\left(\begin{array}{cc}e&b_1\\b_2&e'\end{array}\right)=e$. Then the
map $d\otimes\id_D:C\otimes_{\min}D\to E\otimes_{\min}D$ is a
well-defined contraction.
Let $p_D:E\otimes_{\min}D\to A\otimes_E D$ be the quotient
map. Then
 $$
\bigl(p_D\circ (d\otimes\id_D)\circ (s_t\otimes_\lambda\id_D)\bigr)x=x
 $$
for any $x\in A\odot D$. The maps $p_D$ and $d\otimes\id_D$ are
contractions and the map $s_t\otimes_\lambda\id_D$ is asymptotically a
contraction, hence the composition of these three maps, $A\otimes_\lambda
D\to A\otimes_E D$ is asymptotically a contraction. But that contradicts
$\|x\|_E >\|x\|_\lambda$.
\end{proof}

\section{Distinctness of the asymptotic norm form the minimal and the
maximal ones}

\begin{prop}
There exist separable $C^*$-algebras $A$ and $D$ such that the norms
$\|\cdot\|_\lambda$ and $\|\cdot\|_{\min}$ on $A\odot D$ are different.

\end{prop}
\begin{proof}
Suppose the contrary, i.e. that for any $A$ and $D$ these two norms
coincide. Take a non-exact separable $C^*$-algebra $D$ and an extension
(\ref{extension}), with $A$ contractible, such that the
sequence $0\to B\otimes_{\min}D\to E\otimes_{\min}D\to A\otimes_{\min}D\to
0$ is not exact. Then by Proposition \ref{exact_seq} the extension is not
semi-invertible, making a contradiction to Theorem \ref{contractible}.
\end{proof}

Since $\|\cdot\|_\lambda\leq\|\cdot\|_\sigma$, the latter norm also
differs from the minimal tensor norm in general.

We know nothing about associativity of the asymptotic tensor norms, but
the following argument shows that either they are not associative or they
look much more like the maximal tensor norm than the minimal one, at least
in some group algebra examples.

Let $G$ be an infinite hyperbolic property $T$ group and let $p\in C^*(G)$
be the projection corresponding to the trivial representation of $G$. Let
$\Delta:{\mathbb C}(G)\to C^*(G)\odot C^*(G)$ denote the diagonal map
given by $g\mapsto g\odot g$, $g\in G$. We use the same notation $\Delta$
for $*$-homomorphisms from $C^*(G)$ to $C^*(G)\otimes_{\max}C^*(G)$
and its quotients (like $C^*_r(G)\otimes_{\min}C^*(G)$) as well.

\begin{cor}
Suppose that the asymptotic tensor product is associative with respect to
commutative $C^*$-algebras, i.e.
$C_0(X)\otimes(A\otimes_\lambda D)=(C_0(X)\otimes A)\otimes_\lambda D$.
Then $\Delta(p)$ in $C^*_r(G)\otimes_\lambda C^*(G)$ is not zero.

\end{cor}
\begin{proof}
Put $A=C^*_r(G)$, $E=C^*(G)$ and consider the extension $0\to B\to
E\to A\to 0$, where $B$ is the kernel of the natural
surjection of the full group $C^*$-algebra onto the reduced one. After
tensoring this extension by $C_0[0,1)$ one gets a semi-invertible (even an
asymptotically split) extension, hence by Proposition \ref{exact_seq} one
has $\|\cdot\|_\lambda\geq\|\cdot\|_{C_0[0,1)\otimes E}$ on
$C_0[0,1)\otimes C^*_r(G)\odot C^*(G)$. Suppose that $\Delta(p)=0$ in
$C^*_r(G)\otimes_\lambda C^*(G)$. Then $\|f\otimes\Delta(p)\|_\lambda=0$
for all $f\in C_0[0,1)$. But the quotient tensor norm is associative, so
$\|f\otimes\Delta(p)\|_{C_0[0,1)\otimes
E}=\|f\|\cdot\|\Delta(p)\|_E=\|f\|$, since $\Delta(p)\neq 0$ in
$C^*_r(G)\otimes_E C^*(G)$ (cf. Lemma 6.2.9 of \cite{Guentner-Higson}).
\end{proof}

Similarly, associativity of the symmetric asymptotic tensor product
implies that $\Delta(p)$ is not zero in $C^*_r(G)\otimes_\sigma C^*_r(G)$.

Recall that $\Delta(p)$ is not zero in $C^*_r(G)\otimes_{\max}C^*_{(r)}(G)$
and it is zero in $C^*_r(G)\otimes_{\min}C^*_{(r)}(G)$ \cite{Guentner-Higson}.

\begin{prop}
Let $D$ be a $C^*$-algebra without the weak expectation property of Lance
\cite{Lance} and let $F_\infty$ denote a free froup on an infinite set of
generators. Then
$C^*(F_\infty)\otimes_{\max}D\neq C^*(F_\infty)\otimes_\lambda D=
C^*(F_\infty)\otimes_{\min}D$.

\end{prop}
\begin{proof}
Recall that $D$ has the weak expectation property if
$D\otimes_{\max}E\subset C\otimes_{\max}E$ canonically for any
$C^*$-algebra $E$ and for any $C^*$-algebra $C$ containing $D$ as
$C^*$-subalgebra.
By Proposition 1.1 of \cite{Kirchberg} one has
$C^*(F_\infty)\otimes_{\max}D\neq C^*(F_\infty)\otimes_{\min}D$, so we
have to check that the asymptotic norm coincides with the minimal one.
Since both norms are cross-norms, it is sufficient to check that they
coincide on a dense set, e.g. on finite sums $c=\sum_i\gamma_i\odot
d_i\in \mathbb C[F_\infty]\odot D$, where all $\gamma_i\in F_\infty$
are words on a finite number of generators, $g_1,\ldots,g_n$, of
$F_\infty$. For any asymptotic representation
$\rho_t:F_\infty\to L(H)$ without loss of generality we can assume that
for big enough $t$ the operators $\rho_t(g_k)$, $k=1,\ldots,n$, are 
unitaries. Since $F_\infty$ is free, this means that in calculation of 
the asymptotic norm for $c$ we can use only genuine representations
$\pi_t$ 
of $F_\infty$ given by
$\pi_t(g_k)=\left\lbrace\begin{array}{cl}\rho_t(g_k),&
{\rm if\ }k\leq n,\\1,&{\rm if\ }k>n
\end{array}\right.$ instead of asymptotic ones, hence the asymptotic norm 
coincides with the minimal one.
\end{proof}

As an example of a $C^*$-algebra without the weak expectation property one 
can use the reduced group $C^*$-algebra of $SL_2(\mathbb Z)$ 
\cite{Kirchberg}. 

Thus we have established the following result.
\begin{thm}
The tensor norm $\|\cdot\|_\lambda$ differs both from the minimal and the 
maximal tensor norms. The tensor norm $\|\cdot\|_\sigma$ differs from the 
minimal tensor norm.

\end{thm}

We don't know if the symmetric asymptotic tensor norm differs from the
maximal one, although their coincidence seems very unlikely. The similar
argument does not work, since it reduces to a long-standing problem if
$C^*(F_\infty)\otimes_{\max}C^*(F_\infty)\neq
C^*(F_\infty)\otimes_{\min}C^*(F_\infty)$ \cite{Kirchberg}.


\vspace{2cm}
\parbox{7cm}{V. M. Manuilov\\
Dept. of Mech. and Math.,\\
Moscow State University,\\
Moscow, 119992, Russia\\
e-mail: manuilov@mech.math.msu.su
}
\hfill
\parbox{6cm}{K. Thomsen\\
Institut for matematiske fag,\\
Ny Munkegade, 8000 Aarhus C,\\
Denmark\\
e-mail: matkt@imf.au.dk
}


\begin{thebibliography}{99}

\bibitem{Anderson}
{J. Anderson.} {\it A $C^*$-algebra $A$ for which $\Ext(A)$ is
not a group}, Ann. Math. {\bf 107} (1978), 455--458.

\bibitem{BDF}
{L. G. Brown, R. G. Douglas, P. A. Fillmore.} {\it Extensions
of $C^*$-algebras and $K$-homology}, Ann.  Math. {\bf 105} (1977),
265--324.

\bibitem{CH}
{A. Connes, N. Higson.}
{\it D\'eformations, morphismes asymptotiques et
$K$-th\'eorie bivariante}, C. R. Acad. Sci. Paris S\'er. I Math. {\bf 311}
(1990), 101--106.

\bibitem{Guentner-Higson}
{E. Guentner, N. Higson.}{\it Group $C^*$-algebras and $K$-theory},
Preprint.

\bibitem{Haagerup}
{U. Haagerup.} {\it Random matrices, free probability and the invariant
subspace problem relative to a von Neumann algebra},
Proc. Internat. Congress of Mathem., Beijing 2002,
v. 1, 273--290.

\bibitem{Kirchberg}
{E. Kirchberg.} {\it On semi-split extensions, tensor products and
exactness of $C^*$-algebras}, Invent. Math. {\bf 112} (1993), 449--489.

\bibitem{Lance}
{L. C. Lance.} {\it On nuclear $C^*$-algebras}, J. Funct. Anal. {\bf 12}
(1973), 157--176.

\bibitem{Loring}
{T. Loring.} {\it Almost multiplicative maps between
$C^*$-algebras}, Operator Algebras and Quantum Field Theory, Rome, 1996,
111--122.

\bibitem{MT2}
{V. M. Manuilov, K. Thomsen.} {\it  Asymptotically split extensions
and E-theory},  Algebra i Analiz 12 (2000) No 5, 142--157 (in Russian);
English translation: St. Petersburg Math. J. {\bf 12} (2001), 819--830.

\bibitem{MT3}
{V. M. Manuilov, K. Thomsen.} {\it The Connes--Higson construction is an
isomorphism}, Preprint.

\bibitem{Wass2}
{S. Wassermann.} {\it Liftings in $C^*$-algebras: a counterexample},
Bull. London Math. Soc. {\bf 9} (1976), 201--202.

\bibitem{Wass1}
{S. Wassermann.} {\it Tensor products of free group $C^*$-algebras},
Bull. London Math. Soc. {\bf 22} (1990), 375--380.

\bibitem{Wass-book}
{S. Wassermann.} {\it Exact $C^*$-algebras and related topics}, Lecture
Notes Series {\bf 19}, Seoul National Univ., 1994.


\end{thebibliography}
\end{document}